\documentclass[12pt,leqno]{article}
\usepackage{amssymb,amsbsy,amsmath,amsfonts,amssymb,amscd}

\usepackage[english]{babel}
\usepackage[T1]{fontenc}
\usepackage{indentfirst}

\newtheorem{statement}{}
\newtheorem{theoreme}[statement]{Theorem}
\newtheorem{lemme}[statement]{Lemma}

\newtheorem{proposition}[statement]{Proposition}
\newtheorem{definition}[statement]{Definition}
\newtheorem{corollaire}[statement]{Corollary}

\newcommand\C{\mathbb C}

\newcommand\D{\mathbb D}

\renewcommand\P{\mathbb P}
\newcommand\eps{\varepsilon}
\newcommand\ind{{\rm 1\kern-.30em I}}

\title{A criterion of weak compactness for operators on subspaces of Orlicz spaces}
\author{Pascal Lef\`evre, Daniel Li,\\ Herv\'e Queff\'elec, Luis Rodr{\'\i}guez-Piazza}

\date{\footnotesize \today}

\begin{document}

\maketitle

\noindent{\bf Abstact.} \emph{We give a criterion of weak compactness for the operators on the Morse-Transue 
space $M^\Psi$, the subspace of the Orlicz space $L^\Psi$ generated by $L^\infty$.}
\bigskip


\section{Introduction and Notation.}

In 1975, C. Niculescu established a characterization of weakly compact operators $T$ from ${\cal C}(S)$, where 
$S$ is a compact space, into a Banach space $Z$ (\cite{N1, N2}, see~\cite{DJT} Theorem 15.2 too): 
$T\colon {\cal C}(S) \to Z$ is weakly compact if and only if there exists a Borel probability measure $\mu$ on 
$S$ such that, for every $\epsilon >0$, there exists a constant 
$C(\epsilon) >0$ such that:
\begin{displaymath}
\qquad \| Tf \| \leq C(\epsilon)\,\|f\|_{L^1(\mu)} + \epsilon\,\|f\|_\infty\,, \quad \forall f\in {\cal C}(S).
\end{displaymath} 
The same kind of result was proved by H. Jarchow for $\C^\ast$-algebras in~\cite{Ja}, and by the first author 
for $A(\D)$ and $H^\infty$ (see~\cite{Pascal}). The criterion for $H^\infty$ played a key role to give an 
elementary proof of the equivalence between weak compactness and compactness for composition operators on 
$H^\infty$. 
\smallskip

Beside these spaces, one natural class of Banach spaces is the class of Orlicz spaces $L^\Psi$. Unfortunately, 
we shall see that the above criterion is in general not true for Orlicz spaces. However, it remains true when we 
restrict ourselves to subspaces of the Morse-Transue space $M^\Psi$. This space is the closure of $L^\infty$ in 
the Orlicz space $L^\Psi$.\par\smallskip

In this paper, we first give a characterization of  the operators from a subspace of $M^\Psi$ which fix no copy 
of $c_0$. When the complementary function of $\Psi$ satisfies $\Delta_2$, that gives a criterion of weak 
compactness. If moreover $\Psi$ satisfies a growth condition, that we call $\Delta^0$, the criterion has a more 
usable formulation, analogous to those described above.\par
As in the case of $H^\infty$ (but this is far less elementary), this new version obtained for subspaces of 
Morse-Transue spaces (Theorem \ref{weak compactness}), combined with a study of generalized Carleson 
measures, may be used to prove the equivalence between weak compactness and compactness for composition 
operators on Hardy-Orlicz spaces (see~\cite{LLQR2}), when $\Psi$ satisfies $\Delta^0$.
\par
However, we think that this characterization has an intrinsic interest for Orlicz spaces, and will be useful not 
only for composition operators (see Remark~5 at the end of the paper).\par

\medskip
In this note, we shall consider Orlicz spaces defined on a probability space 
$(\Omega,\P)$, that we shall assume non purely atomic.\par
By an Orlicz function, we shall understand that $\Psi\colon [0, \infty] \to [0,\infty]$ 
is a non-decreasing convex function such that $\Psi(0)=0$ and $\Psi(\infty)=\infty$. To avoid 
pathologies, we shall assume that we work with an Orlicz function $\Psi$ having the 
following additional properties: $\Psi$ is continuous at $0$, strictly convex (hence 
\emph{strictly} increasing), and such that
\begin{displaymath}
\frac{\Psi(x)}{x}\mathop{\longrightarrow}_{x\to \infty} \infty.
\end{displaymath}
This is essentially to exclude the case of $\Psi(x) = ax$. 
The Orlicz space $L^\Psi(\Omega)$ is the space of all (equivalence classes of) measurable
functions $f\colon \Omega\to\C$ for which there is a constant $C > 0$ such that
\begin{displaymath}
\int_\Omega \Psi\Big(\frac{\vert f(t)\vert} {C}\Big)\,d\P(t) < +\infty
\end{displaymath}
and then $\Vert f\Vert_\Psi$ (the \emph{Luxemburg norm}) is the infimum of all possible 
constants $C$ such that this integral is $\leq 1$.\par
To every Orlicz function is associated the complementary Orlicz function 
$\Phi=\Psi^\ast\colon [0,\infty] \to [0,\infty]$ defined by:
\begin{displaymath}
\Phi(x)=\sup_{y\geq 0} \big(xy - \Psi(y)\big).
\end{displaymath}
The extra assumptions on $\Psi$ ensure that $\Phi$ is itself strictly convex.\par
\smallskip

Throughout this paper, we shall assume, except explicit mention of the contrary, that the 
\emph{complementary} Orlicz function satisfies the $\Delta_2$ condition 
($\Phi\in \Delta_2$), \emph{i.e.}, for some constant $K>0$, and some $x_0>0$, we have:
\begin{displaymath}
\Phi(2x)\leq K\,\Phi(x),\hskip 1cm \forall x\geq x_0.
\end{displaymath}
This is usually expressed by saying that $\Psi$ satisfies the $\nabla_2$ condition 
($\Psi\in \nabla_2$). This is 
equivalent to say that for some $\beta>1$ and $x_0>0$, one has 
$\Psi(x)\leq \Psi(\beta x)/(2\beta)$ for $x\geq x_0$, and that implies that 
$\frac{\Psi(x)}{x}\mathop{\longrightarrow}\limits_{x\to \infty} \infty$. In particular, 
this excludes the case $L^\Psi = L^1$.\par\smallskip

When $\Phi$ satisfies the $\Delta_2$ condition, $L^\Psi$ is the dual space of $L^\Phi$.
\par\medskip

We shall denote by $M^\Psi$ the closure of $L^\infty$ in $L^\Psi$. Equivalently (see~\cite{Rao}, 
page 75), $M^\Psi$ is the space of (classes of) functions such that:
\begin{displaymath}
\int_\Omega \Psi\Big(\frac{\vert f(t)\vert} {C}\Big)\,d\P(t)<+\infty,\hskip 1cm 
\forall C>0.
\end{displaymath}
This space is the \emph{Morse-Transue space} associated to $\Psi$, and  
$(M^\Psi)^\ast = L^\Phi$, isometrically if $L^\Phi$ is provided with the Orlicz norm, and 
isomorphically if it is equipped with the Luxemburg norm (see~\cite{Rao}, Chapter IV, 
Theorem 1.7, page 110).\par
We have $M^\Psi = L^\Psi$ if and only if $\Psi$ satisfies the $\Delta_2$ condition, and 
$L^\Psi$ is reflexive if and only if both $\Psi$ and $\Phi$ satisfy the $\Delta_2$ condition. 
When the complementary function $\Phi=\Psi^\ast$ of $\Psi$ satisfies it  
(but $\Psi$ does not satisfy this $\Delta_2$ condition, to exclude the reflexive case), we have 
(see~\cite{Rao}, Chapter IV, Proposition 2.8, page 122, and Theorem 2.11, page 123):
\begin{equation}\label{M-ideal}
(L^\Psi)^\ast= (M^\Psi)^\ast \oplus_1 (M^\Psi)^\perp,\tag*{($\ast$)}
\end{equation}
or, equivalently, $(L^\Psi)^\ast= L^\Phi \oplus_1 (M^\Psi)^\perp$, isometrically, with the 
Orlicz norm on $L^\Phi$.
\par\bigskip

For all the matter about Orlicz functions and Orlicz spaces, we refer to~\cite{Rao}, or to~\cite{Kras}.

\medskip

\noindent{\bf Acknowledgement.} This work was made during the stay in Lens as 
\emph{Professeur invit\'e de l'Universit\'e d'Artois}, in May--June 2005, and a visit in Lens and Lille 
in March 2006 of the fourth-named author. He would like to thank both Mathematics Departments for their 
kind hospitality.
\medskip\goodbreak

\section{Main result.}

Our goal in this section is the following criterion of weak compactness for operators. We begin with:

\begin{theoreme}\label{fixe pas c0}
Let $\Psi$ be an arbitrary Orlicz function, and let $X$ be a subspace of the Morse-Transue space $M^\Psi$. 
Then an operator $T\colon X \to Y$ from $X$ into a Banach space $Y$ fixes no copy of $c_0$ if 
and only if :
\begin{equation}\label{critere c0}
\forall \eps >0 \hskip 3mm \exists C_\eps >0\,: \hskip 5mm 
\|Tf\| \leq \bigg[C_\eps  \int_\Omega \Psi\Big(\eps \frac{|f|}{\ \|f \|_\Psi}\Big)\,d\P 
+ \eps \bigg] \,\|f\|_\Psi,\hskip 2mm \forall f\in X.
\end{equation}
\end{theoreme}

Recall that saying that $T$ fixes a copy of $c_0$ means that there exists a subspace $X_0$ of $X$ isomorphic 
to $c_0$ such that $T$ realizes an isomorphism between $X_0$ and $T(X_0)$.
\medskip

Before proving that, we shall give some consequences. First, we have:

\begin{corollaire}\label{corollaire}
Assume that the complementary function of  $\Psi$ has $\Delta_2$ ($\Psi\in \nabla_2$). Then for 
every subspace $X$ of $M^\Psi$, and every operator $T \colon X \to Y$, $T$ is weakly compact if and 
only if it satisfies \eqref{critere c0}.
\end{corollaire}

\noindent{\bf Proof.} When the complementary function of $\Psi$ has $\Delta_2$, one has the 
decomposition \ref{M-ideal}, which means that $M^\Psi$ is $M$-ideal in its bidual 
(see~\cite{HWW}, Chapter III; this result was first shown by D. Werner (\cite{Werner} -- see 
also~\cite{HWW}, Chapter III, Example 1.4 (d), page 105 -- by a different way, using the ball intersection 
property; note that in these references, it is moreover assumed that $\Psi$ does not satisfy the 
$\Delta_2$ condition, but if it satisfies it, the space $L^\Psi$ is reflexive, and so the result is obvious). 
But every subspace $X$ of a Banach space which is $M$-ideal of its bidual has Pe{\l}czy\'nki's property $(V)$ 
(\cite{GS1, GS2}; see also~\cite{HWW}, Chapter III, Theorem 3.4), which means that operators from $X$ 
are weakly compact if and only if they fix no copy of $c_0$.\hfill$\square$
\medskip

With $\Psi$ satisfying the following growth condition, the characterization \eqref{critere c0} takes on a 
more usable form.

\begin{definition}\label{Delta0}
We say that the Orlicz function $\Psi$ satisfies the \emph{$\Delta^0$ condition} if  for some 
$\beta >1$, one has: 
\begin{displaymath}
\lim_{x\to +\infty} \frac{\Psi(\beta x)}{\Psi(x)}=+\infty.
\end{displaymath}
\end{definition}

\noindent This growth condition is a strong negation of the $\Delta_2$ condition and it implies 
that the complementary function $\Phi=\Psi^\ast$ of $\Psi$ satisfies the $\Delta_2$ condition.\par
Note that in the following theorem, we cannot content ourselves with $\Psi\notin \Delta_2$
(\emph{i.e.} $\limsup_{x \to +\infty} \Psi (\beta x)/\Psi (x) =+\infty$), instead of $\Psi \in \Delta^0$ 
(see Remark~3 in Section~3). An interesting question is whether the condition $\Psi \in \Delta^0$ is 
actually necessary for this characteriztion.

\begin{theoreme}\label{weak compactness} 
Assume that $\Psi$ satisfies the $\Delta^0$ condition, and let $X$ be a subspace of $M^\Psi$. 
Then every linear operator $T$ mapping $X$ into some Banach space $Y$ is weakly compact if and only if 
for some (and then for all) $1\leq p < \infty$:  
\begin{equation}\label{critere}
\forall\eps>0,\ \exists C_\eps>0,\hskip 5mm \|T(f)\|\le C_\eps\|f\|_p+\eps\,\|f\|_\Psi,
\hskip 5mm \forall f\in X.\tag*{(W)}
\end{equation}
\end{theoreme}

\noindent{\bf Remark 1.} This theorem extends~\cite{LLQR} Theorem II.1. As in the case of $\C^\ast$-algebras 
(see~\cite{DJT}, Notes and Remarks, Chap. 15), there are miscellaneous applications of such a characterization.
\medskip

\noindent{\bf Remark 2.} Contrary to the $\Delta_2$ condition where the constant $2$ may be replaced 
by any constant $\beta >1$, in this $\Delta^0$ condition, the constant $\beta$ cannot be replaced by 
another, as the following example shows.
\medskip

\noindent{\bf Example.} There exists an Orlicz function $\Psi$ such that:
\begin{equation}\label{beta=5}
\lim_{x\to +\infty} \frac{\Psi(5x)}{\Psi(x)} = +\infty ,
\end{equation}
but 
\begin{equation}\label{beta not 2}
\liminf_{x\to +\infty} \frac{\Psi(2x)}{\Psi(x)} < +\infty. 
\end{equation}

\noindent{\bf Proof.} Let $(c_n)_n$ be an increasing sequence of positive 
numbers such that  $\displaystyle \lim_{n\to\infty} \frac{c_{n +1}} {c_n}=+\infty$, take 
$\psi(t)= c_n$ for $t\in (4^n,4^{n+1}]$ and $\Psi(x) =\int_0^x \psi(t)\,dt$. Then 
\eqref{beta=5} is verified. On the other hand, if $x_n=2\cdot 4^n$, one has 
$\Psi(x_n)\geq c_n 4^n$, and  $\Psi(2x_n)\leq  c_n 4^{n+1}$, so we get \eqref{beta not 2}.
\hfill$\square$
\medskip

Before proving Theorem~\ref{weak compactness}, let us note that it has the following straightforward corollary.

\begin{corollaire}\label{famille}
Let $X$ be like in Theorem~\ref{weak compactness}, and assume  that ${\cal F}$ is a family of operators from 
$X$ into a Banach space $Y$ with the following property\/: there exists a bounded sequence $(g_n)_n$ in $X$ 
such that $\lim_{n\to\infty} \|g_n\|_1=0$ and such that an operator $T \in {\cal F}$ is compact whenever 
\begin{displaymath}
\lim_{n\to\infty} \| T g_n\|=0.
\end{displaymath}
Then every  weakly compact operator in $T\in {\cal F}$ is  actually compact.
\end{corollaire}

In the forthcoming paper~\cite{LLQR2}, we prove, using a generalization of the notion of Carleson measure, that 
a composition operator $C_\phi \colon H^\Psi \to H^\Psi$ ($H^\Psi$ is the space of analytic functions on 
the unit disk $\D$ of the complex plane whose boundary values are in $L^\Psi (\partial \D)$, and 
$\phi\colon \D \to \D$ is an analytic self-map) is compact whenever:
\begin{displaymath}
\lim_{\ \,r\to 1^-} \sup_{|\xi|=1} \Psi^{-1} \big(1/(1-r)\big) \|C_\phi (u_{\xi,r})\|_\Psi =  0,
\end{displaymath}
where:
\begin{displaymath}
u_{\xi,r}(z) =\Big(\frac{1-r}{1 - \bar{\xi} r z}\Big)^2, \hskip 3mm |z|<1,
\end{displaymath}
and we have:
\begin{displaymath}
 \lim_{\ \,r\to 1^-} \sup_{|\xi|=1} \Psi^{-1} \big(1/(1-r)\big) \| u_{\xi,r}\|_1 =  0
\end{displaymath}
when $C_\phi$ is weakly compact and $\Psi \in \Delta^0$.\\
Though the situation does not fit exactly as in Corollary~\ref{famille} (not because of the space $H^\Psi$, which 
is not a subspace of $M^\Psi$: we actually work in $HM^\Psi = H^\Psi \cap M^\Psi$ since 
$u_{\xi,r} \in HM^\Psi$, but because of the fact that we ask a uniform limit for $|\xi|=1$), the same ideas 
allow us to get, when $\Psi$ satisfies the condition $\Delta^0$, that $C_\phi$ is compact if and only if it is 
weakly compact.
\medskip\goodbreak

\noindent{\bf Proof of Theorem~\ref{weak compactness}.} Assume that we have \ref{critere}. We may 
assume that $p>1$, since if \ref{critere} is satisfied for some $p\geq 1$, it is satisfied for all $q\geq p$. 
Moreover, we may assume that $L^\Psi \stackrel{j}{\hookrightarrow} L^p$ since $\Psi$ satisfies condition 
$\Delta^0$ (since we have: $\lim_{x\to + \infty} \frac{\Psi(x)}{x^r}= +\infty$, for every $r>0$). Then 
$T \big[(1/C_\eps) B_{L^p} \cap (1/\eps) B_X\big] \subseteq 2 B_Y$. Taking the polar of these sets, 
we get $T^\ast (B_{Y^\ast}) \subseteq (2 C_\eps) B_{j^\ast[(L^p)^\ast]} + (2\eps) B_{X^\ast}$, for 
every $\eps >0$. By a well-known lemma of Grothendieck, 
we get, since $B_{j^\ast[(L^p)^\ast]}$ is weakly compact, that $T^\ast(B_{Y^\ast})$ is relatively 
weakly compact, {\it i.e.} $T^\ast$, and hence also $T$, is weakly compact.
\par

Conversely, assume in Theorem \ref{weak compactness} that $T$ is weakly compact. We are going to show 
that \ref{critere} is satisfied with $p=1$ (hence for all finite $p\geq 1$). Let $\eps >0$. Since the 
$\Delta^0$ condition implies that the complementary function of $\Psi$ satisfies the $\Delta_2$ 
condition, Corollary \ref{corollaire} implies that, when $\|f\|_\Psi=1$:
\begin{displaymath}
\|Tf\| \leq C_{\eps/2} \int_\Omega \Psi \big( (\eps/2) |f|\big)\,d\P +\eps/2.
\end{displaymath}
As $\Psi$ satisfies the $\Delta^0$ condition, there is some $\beta >1$ such that 
$\displaystyle \frac{\Psi(x)}{\Psi(\beta x)} \to 0$ as $x\to \infty$; hence, with 
$\kappa= \eps/2C_{\eps/2}$, there exists some $x_\kappa >0$ such that 
$\Psi(x)\leq \kappa \Psi (\beta x)$ for $x\geq x_\kappa$. By the convexity of $\Psi$, one has 
$\Psi(x) \leq \frac{\Psi(x_\kappa)}{x_\kappa} x =: K_\kappa x$ for $0\leq x\leq x_\kappa$. Hence, 
for every $x \geq 0$:
\begin{displaymath}
\Psi(x) \leq \kappa \Psi(\beta x) + K_\kappa x.
\end{displaymath}
It follows that, for $f\in X$, with $\|f\|_\Psi=1$:
\begin{displaymath}
\int_\Omega \Psi \big((\eps/2) |f|\big)\,d\P 
\leq \kappa \int_\Omega \Psi\big(\beta (\eps/2) |f|\big)\,d\P + K_\kappa \frac{\eps}{2} \|f\|_1 
\leq \kappa + K_\kappa \frac{\eps}{2} \|f\|_1
\end{displaymath}
if we have chosen $\eps \leq 2/\beta$. Hence:
\begin{displaymath}
\| Tf\|  \leq C_{\eps/2} \Big( \kappa + K_\kappa \frac{\eps}{2} \|f\|_1\Big) +\frac{\eps}{2} 
= C_{\eps/2} \frac{\eps}{2} K _\kappa \|f\|_1 + \Big(C_{\eps/2} \kappa +\frac{\eps}{2}\Big) 
= C'_{\eps} \| f\|_1 + \eps,
\end{displaymath}
which is \ref{critere}.\hfill$\square$
\medskip

\noindent{\bf Remark.} The sufficient condition is actually a general fact, 
which is surely well known (see~\cite{Pascal}, Theorem 1.1, for a similar result, 
and~\cite{DJT}, Theorem 15.2 for ${\cal C}(K)$; see also~\cite{JM}, page 81), and has close connection 
with interpolation (see~\cite{Beauzamy}, Proposition 1), but we have found no reference, and so we shall 
state it separately without proof (the proof  follows that given in~\cite{DJT}, page 310).

\begin{proposition}\label{general}
Let $T\colon X \to Y$ be an operator between two Banach spaces. Assume that there is a Banach 
space $Z$ and a weakly compact map $j\colon X \to Z$ such that: for every 
$\eps>0$, there exists $C_\eps>0$ such that
\begin{displaymath}
\hskip 1cm 
\Vert Tx \Vert \leq C_\eps \Vert jx\Vert_Z + \eps\,\Vert x \Vert_X\,,\hskip 3mm \forall x\in X.
\end{displaymath}
Then $T$ is weakly compact.
\end{proposition}

Note that, by the Davis-Figiel-Johnson-Pe{\l}czy\'nski factorization theorem, we may assume that 
$Z$ is reflexive. We may also assume that $j$ is injective, because 
$\ker j \subseteq \ker T$, so $T$ induces a map 
$\tilde T\colon X/\ker j \to Y$ with the same property as $T$. Indeed, if $jx=0$, then 
$\|Tx\|\leq \eps \|x\|$ for every $\eps >0$, and hence $Tx=0$.\par
\bigskip

\noindent{\bf Proof of Theorem \ref{fixe pas c0}.}  Assume first that $T$ fixes a copy of $c_0$. There are  
hence  some $\delta >0$ and a sequence $(f_n)_n$ in $X$ equivalent to the canonical basis of $c_0$ 
such that $\|f_n\|_\Psi=1$ and $\|Tf_n\|\geq \delta$. In particular, there is some $M>0$ such that, for 
every choice of $\eps_n =\pm 1$:
\begin{displaymath}
\Big\|\sum_{n=1}^N \eps_n f_n \Big\|_\Psi \leq M, \hskip 5mm \forall N\geq 1.
\end{displaymath}
Let $(r_n)_n$ be a Rademacher sequence. We have, first by Khintchine's inequality, then by 
Jensen's inequality and Fubini's Theorem:
\begin{align*}
\int_\Omega \Psi\bigg(\frac{1}{M\sqrt{2}} \Big(\sum_{n=1}^N |f_n|^2\Big)^{1/2} \bigg)\,d\P
& \leq \int_\Omega \Psi\bigg[ \frac{1}{M} \int_0^1 \Big|\sum_{n=1}^N r_n(t) f_n \Big|\,dt\bigg]\,d\P \\
& \leq \int_\Omega\int_0^1 \Psi\bigg[ \frac{1}{M} \Big|\sum_{n=1}^N r_n(t) f_n \Big|\,dt\bigg]\,d\P \\
& = \int_0^1\int_\Omega \Psi\bigg[ \frac{1}{M} \Big|\sum_{n=1}^N r_n(t) f_n \Big|\,d\P\bigg]\,dt 
\leq 1.
\end{align*}
The monotone convergence Theorem gives then:
\begin{displaymath}
\int_\Omega \Psi\bigg(\frac{1}{M\sqrt{2}} \Big(\sum_{n=1}^\infty |f_n|^2\Big)^{1/2} \bigg)\,d\P 
\leq 1.
\end{displaymath}
In particular, $\sum_{n=1}^\infty |f_n|^2$ is finite almost everywhere, and hence $f_n \to 0$ 
almost everywhere. Since 
$ \Psi\Big(\frac{1}{M\sqrt{2}} \big(\sum_{n=1}^\infty |f_n|^2\big)^{1/2} \Big) \in L^1$, by the 
above inequalities, Lebesgue's dominated convergence Theorem gives: 
\begin{displaymath}
\int_\Omega \Psi\Big(\frac{|f_n|}{M\sqrt{2}}\Big)\,d\P \mathop{\longrightarrow}_{n\to\infty} 0.
\end{displaymath}
But that contradicts \eqref{critere c0} with $\eps \leq 1/M\sqrt{2}$ and $\eps < \delta$, since 
$\|Tf_n\| \geq \delta$.
\medskip

The converse follows from the following lemma.

\begin{lemme}\label{construction}
Let $X$ be a subspace of $M^\Psi$, and let $(h_n)_n$ be a sequence in $X$, with 
$\| h_n\|_\Psi=1$ for all $n\geq 1$, and such that, for some $M >0$:
\begin{displaymath}
\int_\Omega \Psi (|h_n|/M)\,d\P \mathop{\longrightarrow}\limits_{n\to \infty} 0.
\end{displaymath}
Then $(h_n)_n$ has a subsequence equivalent to the canonical basis of $c_0$.
\end{lemme}

Indeed, if  condition \eqref{critere c0} is not satisfied, there exist some $\eps_0>0$ and functions 
$h_n \in X$ with $\|h_n\|_\Psi=1$ such that 
$\|Th_n\| \geq 2^n \int_\Omega \Psi (\eps_0 |h_n|)\,d\P +\eps_0$. That implies that 
$\int_\Omega \Psi (\eps_0 |h_n|)\,d\P$ tends to $0$, so Lemma \ref{construction} ensures that 
$(h_n)_n$ has a subsequence, which we shall continue to denote by  $(h_n)_n$, equivalent to the 
canonical basis of $c_0$. Then $(Th_n)_n$ is weakly unconditionally Cauchy. Since 
$\|Th_n\|  \geq \eps_0$, $(Th_n)_n$ has, by Bessaga-Pe{\l}czy\'nski's Theorem, a further subsequence 
equivalent to the canonical basis of $c_0$. It is then obvious that $T$ realizes an isomorphism 
between the spaces generated by these subsequences.\hfill$\square$
\medskip

\noindent{\bf Proof of Lemma \ref{construction}.} The proof uses the idea of the construction made in the 
proof of Theorem II.1 in~\cite{LLQR}, which it generalizes, but with some additional details.\par
By the continuity of $\Psi$, there exists $a>0$ such that $\Psi(a)=1$. Then, since $\Psi$ is 
increasing, we have, for every $g\in L^\infty$:
\begin{displaymath}
\int_\Omega\Psi\Big(a\,\frac{\vert g\vert\hskip 3mm}{\Vert g\Vert_\infty}\Big)
\,d\P\leq 1\,,
\end{displaymath}
and so $\Vert g\Vert_\Psi\leq (1/a)\,\Vert g\Vert_\infty$.\par

Now, choose, for every $n\geq 1$, positive numbers $\alpha_n < a/2^{n+2}$ such that 
$\Psi (\alpha_n/2M) \leq 1$.\par\smallskip

We are going to construct inductively a subsequence $(f_n)_n$ of $(h_n)_n$, a sequence of functions 
$g_n\in L^\infty$ and two sequences of positive numbers $\beta_n$ and 
$\eps_n\leq \min\{ 1/2^{n+1}, M/2^{n+1}\}$, 
such that, for every $n\geq 1$:\smallskip

\begin{itemize}
\item[(i)] if we set $M_1=1$ and, for $n\geq 2$:
$$M_n =
\max\Big\{1, \Psi\Big(\frac{\Vert g_1\Vert_\infty+\cdots+\Vert g_{n-1}\Vert_\infty}{2M}\Big)\Big\}\,,$$
then $M_n \beta_n\leq 1/2^{n+1}$;
\item[(ii)] $\Vert f_n\Vert_\Psi=1$;
\item[(iii)] $\Vert f_n - g_n\Vert_\Psi \leq \eps_n$, with $\eps_n$ such that 
$\beta_n \Psi(\alpha_n/2\eps_n)\geq 2$;
\item[(iv)] $\P(\{\vert g_n\vert > \alpha_n\})\leq \beta_n$;
\item[(v)] $\Vert \breve g_n\Vert_\Psi\geq 1/2$, with 
$\breve g_n=g_n\,\ind_{\{\vert g_n\vert>\alpha_n\}}$.
\end{itemize}
\smallskip

We shall give only the inductive step, since the starting one unfolds identically. 
Suppose hence that the functions $f_1,\ldots,f_{n-1}$, $g_1,\ldots, g_{n-1}$ and the numbers 
$\beta_1,\ldots,\beta_{n-1}$ and $\eps_1,\ldots,\eps_{n-1}$ have been constructed. 
Choose then $\beta_n>0$ such  that $M_n \beta_n\leq 1/2^{n+1}$. Note that $M_n \geq 1$ implies that 
$\beta_n \leq 1/2^{n+1}$.\par

Since $\int_\Omega \Psi(|h_k|/M)\,d\P \to 0$ as $n\to \infty$, we can find $f_n = h_{k_n}$ 
such that $\Vert f_n\Vert_\Psi=1$, and  moreover:
\begin{displaymath}
\P(\{\vert f_n\vert>\alpha_n/2\})
\leq \frac{1}{\Psi(\alpha_n/2M)} \int_\Omega \Psi \Big(\frac{ |f_n|}{M}\Big)\,d\P 
\leq \frac{\beta_n}{2}\cdot
\end{displaymath}
Take now $\eps_n\leq \min\{1/2^{n+1}, M/2^{n+1}\}$ such that 
$0 < \eps_n \leq \alpha_n/ 2\Psi^{-1}(2/\beta_n)$ and $g_n\in L^\infty$ such that 
$\Vert f_n - g_n \Vert_\Psi\leq \eps_n$. Then, since 
\begin{displaymath}
\P(\{|f_n -g_n| > \alpha_n/2\}) \Psi\Big( \frac{\alpha_n}{2\eps_n}\Big) \leq  
\int_\Omega \Psi\Big(\frac{|f_n - g_n|}{\eps_n}\Big)\,d\P \leq 1,
\end{displaymath}
we have:
\begin{align*}
\P(\{|g_n| > \alpha_n\}) &\leq \P(\{|f_n| > \alpha_n/2\}) + \P(\{|f_n -g_n| > \alpha_n/2\}) \\
& \leq \frac{\beta_n}{2} + \frac{1}{\Psi(\alpha_n/2\eps_n)} \leq \beta_n.
\end{align*}
To end the construction, it remains to note that
\begin{align*}
\Vert f_n - \breve g_n\Vert_\Psi &\leq \Vert f_n - g_n\Vert_\Psi + 
\Vert \breve g_n - g_n\Vert_\Psi 
\leq \eps_n + \frac{1}{a} \Vert \breve g_n - g_n\Vert_\infty \\ 
&\leq \frac{1}{2^{n+1}} + \frac {\alpha_n}{a} \leq \frac{1}{2^n}\leq \frac{1}{2}
\end{align*}
and so:
\begin{displaymath}
\Vert \breve g_n\Vert_\Psi\geq \Vert f_n\Vert_\Psi - \Vert f_n - \breve g_n\Vert_\Psi 
\geq 1 - \frac{1}{2} = \frac{1}{2}\cdot
\end{displaymath}
This ends the inductive construction.
\smallskip

Consider now 
\begin{displaymath}
\breve g=\sum_{n=1}^{+\infty} \vert \breve g_n\vert\,.
\end{displaymath}
Set $A_n=\{\vert g_n\vert>\alpha_n\}$ and, for $n\geq 1$:
\begin{displaymath}
B_n=A_n\setminus \bigcup_{j>n} A_j\,.
\end{displaymath}
We have $\P\big(\limsup A_n\big)=0$, because 
\begin{displaymath}
\sum_{n\geq 1} \P(A_n)\leq \sum_{n\geq 1} \beta_n \leq \sum_{n\geq 1} \frac{1}{2^n}<+\infty.
\end{displaymath} 
Now $\breve g$ vanishes out of 
$\bigcup\limits_{n\geq 1} B_n\cup \big(\limsup A_n\big)$ and we have:
\begin{align*}
\int_{B_n} \Psi\Big(\frac{\vert \breve g_n\vert}{2M}\Big)\,d\P 
& \leq \int_\Omega \Psi\Big(\frac{\vert g_n\vert}{2M}\Big)\,d\P \\
& \leq \int_\Omega \Psi\Big(\frac{\vert g_n -f_n\vert}{2M} + \frac{|f_n|}{2M}\Big)\,d\P \\
& \leq \frac{1}{2} \int_\Omega \Psi\Big(\frac{\vert g_n -f_n\vert}{M} \Big)\,d\P 
+ \frac{1}{2} \int_\Omega \Psi\Big(\frac{\vert f_n\vert}{M} \Big)\,d\P. 
\end{align*}
The first integral is less than $\eps_n/M$, because $\Psi (at) \leq a\Psi (t)$ for $0\leq a \leq 1$ 
and $\eps_n/M \leq 1$, so that:
\begin{displaymath}
\int_\Omega \Psi\Big(\frac{\vert g_n -f_n\vert}{M} \Big)\,d\P 
\leq \frac{\eps_n}{M}\int_\Omega \Psi\Big(\frac{\vert g_n -f_n\vert}{\eps_n} \Big)\,d\P 
\leq \frac{\eps_n}{M} \le\frac{1}{2^{n+1}}
\end{displaymath}
(since $\| f_n - g_n \|_\Psi \leq \eps_n$).\par
Since:
\begin{displaymath}
\int_\Omega \Psi\Big(\frac{\vert f_n\vert}{M} \Big)\,d\P \leq \frac{\beta_n}{2} \Psi \big( \alpha_n/ 2M) 
\leq \beta_n/2,
\end{displaymath}
we obtain:
\begin{displaymath}
\int_{B_n}  \Psi\Big( \frac{\vert \breve g_n\vert}{2M}\Big)\,d\P 
\leq\frac{1}{2^{n+2}} + \frac{\beta_n}{4}\,\cdot
\end{displaymath}
Therefore, since $\P (B_n) \leq \P (A_n) \leq \beta_n$:
\begin{align*}
\int_\Omega \Psi\Big(\frac{\vert \breve g\vert}{4M}\Big)\,d\P
&=\sum_{n=1}^{+\infty} \int_{B_n} \Psi\Big(\frac{\vert \breve g\vert}{4M}\Big)\,d\P \cr
&\leq \sum_{n=1}^{+\infty} \int_{B_n} \frac{1}{2}
\bigg[\Psi\Big(\frac{\Vert g_1\Vert_\infty+\cdots+\Vert g_{n-1}\Vert_\infty}{2M}\Big)
+\Psi\Big(\frac{\vert \breve g_n\vert}{2M}\Big)\bigg]\,d\P \cr
&\hskip 5mm
\hbox{by convexity of $\Psi$ and because $\breve g_j=0$ on $B_n$ for $j>n$} \cr
& \leq \frac{1}{2}\sum_{n=1}^{+\infty} 
\Big( M_n\beta_n + \frac{1}{2^{n+2}} + \frac{\beta_n}{4} \Big)  \cr
& \leq \frac{1}{2}\sum_{n=1}^{+\infty} 
\Big(\frac{1}{2^{n+1}} + \frac{1}{2^{n+2}} + \frac{1}{2^{n+2}}\Big)\le1.
\end{align*}
That proves that $\breve g\in L^\Psi$, and consequently that the series 
$\sum_{n\geq 1} \breve g_n$ is weakly unconditionally Cauchy in $L^\Psi$:
\begin{displaymath}
\sup_{n\geq 1}\sup_{\theta_k=\pm 1} 
\Big\Vert \sum_{k=1}^n \theta_k \breve g_k \Big\Vert_\Psi \leq 
\sup_{n\geq 1}\Big\Vert \sum_{k=1}^n |\breve g_k|\, \Big\Vert_\Psi \leq 
\Vert \breve g \Vert_\Psi \leq 4M.
\end{displaymath}
\smallskip

Since $\Vert \breve g_n\Vert_\psi\geq 1/2$, $(\breve g_n)_{n\geq 1}$ has, by 
Bessaga-Pe\l czy\'nski's theorem, a subsequence $(\breve g_{n_k})_{k\geq 1}$ which is 
equivalent to the canonical basis of $c_0$. The corresponding subsequence 
$(f_{n_k})_{k\geq 1}$ of $(f_n)_{n\geq 1}$ remains equivalent to the canonical basis of 
$c_0$, since
\begin{displaymath}
\sum_{n=1}^{+\infty} \Vert f_n - \breve g_n\Vert_\Psi \leq 
\sum_{n=1}^{+\infty} \eps_n +\frac{\alpha_n}{a} \leq 
\sum_{n=1}^{+\infty} \frac{1}{2^{n+1}} + \frac{1}{2^{n+2}} < 1\,.
\end{displaymath}
That ends the proof of Lemma \ref{construction}.\hfill$\square$

\section{Comments}

\noindent{\bf Remark 1.} Let us note that the assumption $X\subseteq M^\Psi$ cannot be relaxed in general. 
In fact, suppose that $X$ is a subspace of $L^\Psi$ containing $L^\infty$, and let 
$\xi\in (M^\psi)^\perp \subseteq (L^\Psi)^\ast$. Being of rank one, $\xi$ is trivially 
weakly compact. Suppose that it satisfies~\ref{critere}. Let 
$f\in X$ with norm $1$, and let $\eps>0$. For $t$ large enough and 
$f_t= f\ind_{\{|f|\leq t\}}$, we have $\Vert f - f_t \Vert_2\leq \eps/C_\eps$. 
Moreover, $f_t \in L^\infty\subseteq X$ and $\| f_t \|_\Psi \leq \| f\|_\Psi=1$. 
Since $\xi$ vanishes on $L^\infty$ and $f - f_t\in X$, we get:
\begin{displaymath}
|\xi(f)| =|\xi (f - f_t)| \leq C_\eps \| f - f_t\|_2 +\eps \| f - f_t\|_\Psi \leq 3\eps.
\end{displaymath}
This implies that $\xi(f)=0$. Since this occurs for every $\xi\in (M^\Psi)^\perp$, we get that 
$X\subseteq M^\Psi$ (and actually $X=M^\Psi$ since $X$ contains $L^\infty$).\hfill$\square$\par
\smallskip
In particular Theorem~\ref{weak compactness} does not hold for $X=L^\Psi$.
\medskip

\noindent{\bf Remark 2.} However, condition~\ref{critere} remains true for bi-adjoint operators coming 
from subspaces of $M^\Psi$: if $T\colon X\subseteq M^\Psi \to Y$ satisfies the 
condition~\ref{critere}, then $T^{\ast\ast} \colon X^{\ast\ast} \to Y^{\ast\ast}$ also satisfies it. 
Indeed, for every $\eps>0$, we get an equivalent norm $\||\,.\,\,|\|_\eps$ on $X$ by putting:
\begin{displaymath}
\|| f |\|_\eps = C_\eps \| f \|_2 + \eps \| f \|_\Psi.
\end{displaymath}
Hence if $f\in X^{\ast\ast}$, there exists a net $(f_\alpha)_\alpha$ of elements in $X$, 
with $\|| f_\alpha |\|_\eps \leq \|| f |\|_\eps$ which converges weak-star to $f$. Then 
$(Tf_\alpha)_\alpha$ converges weak-star to $T^{\ast\ast} f$, and:
\begin{align*}
\| T^{\ast\ast} f \| & \leq \liminf_\alpha \| Tf_\alpha \|  
\leq \liminf_\alpha (C_\eps \| f_\alpha \|_2 + \eps \| f_\alpha \|_\Psi) \\
& = \liminf_\alpha \|| f_\alpha |\|_\eps \leq \|| f |\|_\eps 
= C_\eps \| f \|_2 + \eps \| f \|_\Psi.
\end{align*}
Hence, from Proposition~\ref{general} below, for such a $T$, $T^{\ast\ast}$ is weakly compact 
if and only if it satisfies~\ref{critere}. We shall use this fact in the forthcoming paper~\cite{LLQR2}.
\bigskip

\noindent{\bf Remark 3.} In Theorem~\ref{weak compactness}, we cannot only assume that 
$\Psi\notin \Delta_2$, instead of $\Psi\in \Delta^0$, as the following example shows. It also shows 
that in Corollary~\ref{corollaire}, we cannot obtain condition~\ref{critere} instead of 
condition \eqref{critere c0}.\par
\medskip

\noindent{\bf Example.} Let us define: 
\begin{displaymath}
\psi(t)= \left\{
\begin{array}{lll}
\,t & \text{for} & 0\leq t <1, \\
(k!)(k+2)t -k!(k+1)! & \text{for} & k! \leq t \leq (k+1)!,\hskip 3mm k\geq 1,
\end{array}
\right.
\end{displaymath}
($\psi(k!)=(k!)^2$ for every integer $k\geq 1$ and $\psi$ is linear between $k!$ and $(k+1)!$), and
\begin{displaymath}
\Psi(x)=\int_0^x \psi(t)\,dt.
\end{displaymath}
Since $t^2 \leq \psi(t)$ for all $t\geq 0$, one has $x^3/3 \leq \Psi(x)$ for all $x\geq 0$.
Then
\begin{displaymath}
\Psi(2.n!) \geq \int_{n!}^{2.n!} \psi(t)\,dt = n!(n+2)\frac{3}{2}(n!)^2 - (n!)^2 (n+1)! 
= (n!)^3 \Big(\frac{n}{2} +2\Big)\,\raise0,5mm\hbox{,}
\end{displaymath}
whereas
\begin{displaymath}
\Psi(n!) = \int_0^{n!} \psi(t)\,dt \leq (n!)^2 \, n!=(n!)^3\,;
\end{displaymath}
hence
\begin{displaymath}
\frac{\Psi(2.n!)}{\Psi(n!)} \geq \frac{n}{2} + 2,
\end{displaymath}
and so 
\begin{displaymath}
\limsup_{x\to +\infty} \frac{\Psi(2x)}{\Psi(x)} =+\infty,
\end{displaymath}
which means that $\Psi\notin \Delta_2$.\par
On the other hand, for every $\beta >1$:
\begin{displaymath}
\Psi(n!/\beta) \geq \frac{1}{3}\Big(\frac{\,n!}{\beta}\Big)^3 = \frac{(n!)^3}{3\beta^3}\,\raise0,5mm \hbox{,}
\end{displaymath}
so
\begin{displaymath}
\frac{\Psi(n!)}{\Psi(n!/\beta)} \leq  \frac{(n!)^3}{(n!)^3/3\beta^3} = 3\beta^3\,;
\end{displaymath}
hence 
\begin{displaymath}
\liminf_{x\to +\infty} \frac{\Psi(2x)}{\Psi(x)} \leq 3\beta^3\,,
\end{displaymath}
and $\Psi\notin \Delta^0$ (actually, this will follow too from the fact that Theorem~\ref{weak compactness} 
is not valid for this $\Psi$).\par
Moreover, the conjugate function of $\Psi$ satisfies the condition $\Delta_2$. Indeed, since $\psi$ is convex, 
one has $\psi(2u)\geq 2\psi(u)$ for all $u\geq 0$, and hence:
\begin{displaymath}
\Psi(2x)=\int_0^{2x} \psi(t)\,dt =2\int_0^x \psi(2u)\,du \geq 2\int_0^x 2\psi(u)\,du = 4 \Psi(x), 
\end{displaymath}
and as it was seen in the Introduction that means that $\Psi\in \nabla_2$.
\smallskip

Now, we have $x^3/3 \leq \Psi(x)$ for all $x\geq 0$; therefore $\|\,.\,\|_3 \leq 3^{1/3} \|\,.\,\|_\Psi$.
In particular, we have an inclusion map $j\colon M^\Psi \hookrightarrow L^3$, which is, of course, 
weakly compact. Nevertheless, assuming that $\P$ is diffuse, condition~\ref{critere} is not verified by $j$, 
when $\eps<1$. Indeed, as we have seen before, one has $\Psi(n!)\leq (n!)^3$. Hence, if we choose a 
measurable set $A_n$ such that $\P(A_n)=1/\Psi(n!)$, we have:
\begin{displaymath}
\|\ind_{A_n}\|_\Psi =\frac{1}{\Psi^{-1}\big(1/\P(A_n)\big)} = \frac{1}{n!}\,;
\end{displaymath}
whereas:
\begin{displaymath}
\|\ind_{A_n}\|_3 = \P(A_n)^{1/3} = \frac{1}{\Psi(n!)^{1/3}} \geq \frac{1}{n!}
\end{displaymath}
and
\begin{displaymath}
\|\ind_{A_n}\|_2= \P(A_n)^{1/2} \leq \bigg[\frac{3\ }{(n!)^3}\bigg]^{1/2} = 
\frac{\sqrt 3}{(n!)^{3/2}}\,\cdot
\end{displaymath}
If condition~\ref{critere} were true, we should have, for every $n\geq 1$:
\begin{displaymath}
\frac{1}{n!} \leq C_\eps\, \frac{\sqrt 3\ }{(n!)^{3/2}} + \eps\,\frac{1}{n!}\,\raise0,5mm\hbox{,}
\end{displaymath}
that is:
\begin{displaymath}
\sqrt{n!} \leq \sqrt 3\,\frac{C_\eps}{1-\eps}\,\raise0,5mm\hbox{,}
\end{displaymath}
which is of course impossible for $n$ large enough.\hfill$\square$
\medskip

\noindent{\bf Remark 4.} In the case of the whole space $M^\Psi$, we can give a direct proof 
of the necessity in Theorem~\ref{weak compactness}:\par
\smallskip

\noindent Suppose that $T\colon M^\Psi \to X$ is weakly compact. Then 
$T^\ast \colon X^\ast \to L^\Phi= (M^\Psi)^\ast$ is weakly compact, and so the 
set $K=T^\ast(B_{X^\ast})$ is relatively weakly compact.\par
Since $\Phi$ satisfies the $\Delta^0$ condition, it follows from~\cite{Alexo} (Corollary 2.9) 
that $K$ has equi-absolutely continuous norms. Hence, for every $\eps>0$, we can find 
$\delta_\eps>0$ such that:
\begin{displaymath}
m(A)\leq \delta_\eps \hskip 5mm \Rightarrow \hskip 5mm \Vert g \ind_A\Vert_\Phi\leq \eps/2
\,,\hskip 3mm  \forall g\in T^\ast(B_{X^\ast}).
\end{displaymath}
But (the factor $1/2$ appears because we use the Luxemburg norm on the dual, and not the Orlicz 
norm: see~\cite{Rao}, Proposition III.3.4):
\begin{align*}
\sup_{g\in T^\ast(B_{X^\ast})} \Vert g\ind_A\Vert_\Phi 
&\geq \frac{1}{2}\sup_{u\in B_{X^\ast}} \sup_{\Vert f\Vert_\Psi\leq 1} |<f, (T^\ast u)\ind_A>| 
\\
&= \frac{1}{2}\sup_{u\in B_{X^\ast}} \sup_{\Vert f\Vert_\Psi\leq 1} 
\Big\vert\int f (T^\ast u)\ind_A\,dm \Big\vert \\
&=\frac{1}{2}\sup_{u\in B_{X^\ast}} \sup_{\Vert f\Vert_\Psi\leq 1} \vert <T(f\ind_A),u>\vert 
=\frac{1}{2}\sup_{\Vert f\Vert_\Psi\leq 1} \Vert T(f\ind_A)\Vert; 
\end{align*}
so
\begin{displaymath}
m(A)\leq \delta_\eps \hskip 5mm \Rightarrow \hskip 5mm 
\sup_{\Vert f\Vert_\Psi\leq 1} \Vert T(f\ind_A)\Vert\leq \eps.
\end{displaymath}

Now, we have:

\begin{displaymath}
m(|f|\geq \Vert f\Vert_2/\delta_\eps)\leq \frac{\delta_\eps}{\Vert f\Vert_2}\int |f|\,dm 
=\frac{\delta_\eps}{\Vert f\Vert_2} \Vert f\Vert_1 \leq \delta_\eps;
\end{displaymath}
hence, with $A=\{|f|\geq \Vert f\Vert_2/\delta_\eps\}$, we get, for $\Vert f\Vert_\Psi\leq 1$:
\begin{displaymath}
\Vert Tf \Vert 
\leq \Vert T(f\ind_A)\Vert +\Vert T(f\ind_{A^c})\Vert 
\leq \eps + \Vert T\Vert \frac{\Vert f\Vert_2}{\delta_\eps} 
\end{displaymath}
since $|f\ind_{A^c}|\leq \Vert f\Vert_2/\delta_\eps$ implies 
$\Vert f\ind_{A^c}\Vert_\Psi\leq \Vert f\ind_{A^c}\Vert_\infty 
\leq \Vert f\Vert_2/\delta_\eps$.\hfill $\square$
\medskip

\noindent{\bf Remark~5.} Conversely, E. Lavergne (\cite{Emma}) recently uses  our 
Theorem~\ref{weak compactness} to give a proof of the above quoted result of J. Alexopoulos (\cite{Alexo}, 
Corollary~2.9), and uses it to show that, when $\Psi \in \Delta^0$, then the reflexive subspaces of 
$L^\Phi$ (where $\Phi$ is the conjugate of $\Psi$) are closed for the $L^1$-norm.\par

\bigskip
\vbox{\noindent{\it 
P. Lef\`evre and D. Li, Universit\'e d'Artois,
Laboratoire de Math\'ematiques de Lens EA 2462, F\'ed\'eration CNRS Nord-Pas-de-Calais FR 2956, 
Facult\'e des Sciences Jean Perrin,
Rue Jean Souvraz, S.P.\kern 1mm 18,\par\noindent
62\kern 1mm 307 LENS Cedex,
FRANCE \\ 
pascal.lefevre@euler.univ-artois.fr \hskip 3mm -- \hskip 3mm 
daniel.li@euler.univ-artois.fr
\smallskip

\noindent
H. Queff\'elec,
Universit\'e des Sciences et Techniques de Lille, 
Labo\-ratoire Paul Painlev\'e U.M.R. CNRS 8524, 
U.F.R. de Math\'ematiques,\par\noindent
59\kern 1mm 655 VILLENEUVE D'ASCQ Cedex, 
FRANCE \\ 
queff@math.univ-lille1.fr
\smallskip

\noindent
Luis Rodr{\'\i}guez-Piazza, Universidad de Sevilla, Facultad de 
Matem\'aticas, Dpto de An\'alisis Matem\'atico, Apartado de Correos 1160,\par\noindent 
41\kern 1mm 080 SEVILLA, SPAIN \\ 
piazza@us.es\par}
}

\end{document}